\documentclass[11pt]{article}

\usepackage[T1]{fontenc}
\usepackage[utf8]{inputenc}
\usepackage{lmodern}
\usepackage[english]{babel}
\usepackage{amsmath,amssymb,amsthm}

\newtheorem{thm}{Theorem}
\newtheorem{rem}[thm]{Remark}
\newtheorem{cor}[thm]{Corollary}
\newtheorem{defn}[thm]{Definition}
\newtheorem{lem}[thm]{Lemma}
\newtheorem{prop}[thm]{Proposition}

\usepackage{xcolor}
\definecolor{G}{rgb}{0.0, 0.5, 0.0}
\definecolor{B}{rgb}{0.4, 0.2, 0.1}
\definecolor{BL}{rgb}{0.0, 0.3, 1.0}
\definecolor{L}{rgb}{0.6, 0.0, 0.8}
\definecolor{O}{rgb}{0.8, 0.3, 0.0}
\definecolor{R}{rgb}{0.8, 0.0, 0.0}
\definecolor{LINK}{rgb}{0.35, 0.15, 0.1}
\definecolor{URL}{rgb}{0.6, 0.2, 0.7}
\definecolor{CITE}{rgb}{0.05, 0.3, 0.05}

\newcommand{\TC}{\textcolor}

\usepackage{tikz}
\usetikzlibrary{positioning}

\usepackage{graphicx}
\usepackage{listings}
\usepackage{caption}
\usepackage[linkcolor=LINK, urlcolor=URL, citecolor=CITE]{hyperref}
\usepackage[all]{hypcap}
\hypersetup{linktocpage,colorlinks=true,pdfborder={0 0 0}}

\newcommand{\OEIS}[1]{\href{https://oeis.org/#1}{#1}}
\newcommand{\Modtxt}[1]{(\mathrm{mod}\,#1)} 

\usepackage{geometry}
\title{Deterministic Structures in the Coefficient Stopping-Time Dynamics of the $3x+1$ Problem\\[5mm]}

\author{Mike Winkler\\[5mm]Fakult\"at f\"ur Mathematik\\Ruhr-Universit\"at Bochum, Germany\\mike.winkler@ruhr-uni-bochum.de}
\date{Revised July 27, 2026\\[5mm]}

\begin{document}
	\maketitle
	\begin{abstract}
	    The $3x+1$ problem concerns the iteration of the map $T:\mathbb{Z}\to\mathbb{Z}$ defined by $T(x)=x/2$ for even $x$ and $T(x)=(3x+1)/2$ for odd $x$. We study the coefficient stopping time in the sense of Terras. For each order $n$, we characterize the corresponding residue classes modulo $2^{\sigma_n}$ by the admissible positions of the odd iterates. These position vectors form a directed rooted tree under deletion of the final odd position. This description yields a Pascal-type recursion for the number of classes and proves that the recursive generation produces exactly the admissible vectors, each exactly once.
	    
	    The affine iterate formula gives explicit congruence relations for the classes and arithmetic transition rules between related parity vectors. For every fixed $N$, the union of the classes generated up to order $N$ is periodic with period $2^{\sigma_N}$. Its density is given by an exact finite sum, and its largest initial interval of coverage can be computed from one period. These results concern finite coefficient stopping-time structures. They neither prove that every starting value has finite coefficient stopping time nor establish equality between the coefficient stopping time and the classical stopping time.
	\end{abstract}
	


	\section{Introduction}
	\label{sec:1}
	
	The $3x + 1$ \emph{map} $T:\mathbb{Z}\rightarrow\mathbb{Z}$ is defined by
	\[
		T(x)=\left\{\begin{array}{lcr}\;\;\;\;\displaystyle{\frac{x}{2}} & \mbox{if $x$ is even},\\ \\\displaystyle{\frac{3x+1}{2}} & \mbox{if $x$ is odd}.\end{array}\right.
	\]
	Let $T^0(x)=x$ and $T^s(x)=T\left(T^{s-1}(x)\right)$ for $s\in\mathbb{N}$. For each $x\in\mathbb{N}$, the sequence of iterates $C(x)=\left(T^s(x)\right)_{s=0}^{\infty}$ is referred to as the \emph{trajectory} of $x$. For instance, the initial value $x=11$ yields the trajectory
	\[
		C(11)=(11,17,26,13,20,10,5,8,4,2,1,2,1,2,1,\dotsc).
	\]
	
	If a trajectory does not eventually enter a cycle, then it is unbounded. The $3x + 1$ Conjecture asserts that for every $x\in\mathbb{N}$, the sequence $C(x)$ eventually enters the trivial cycle $(2,1,2,1,\dotsc)$.
	
	We study the coefficient stopping-time structure associated with the iteration of the map $T$. Using parity vectors and the corresponding exponential Diophantine equations, we obtain an explicit recursion for the coefficient stopping-time classes and arithmetic transition rules between adjacent classes in the induced tree. This yields a concrete algebraic description of the coefficient stopping-time classes and their local relations.
	
	Terras~\cite{Terras76,Terras78} and Everett~\cite{Everett77} established foundational density results for stopping-time structures. Related probabilistic and dynamical approaches have been developed by Wirsching~\cite{Wirsching98}. See also the survey by Lagarias~\cite{Lagarias85}. The present work focuses on an explicit finite coefficient stopping-time structure encoded by congruence classes. We establish (i) an exact characterization by admissible odd positions, (ii) a complete rooted-tree construction and a Pascal-type counting recursion, (iii) deterministic arithmetic transition rules between related classes, and (iv) periodicity and exact density formulas for finite-order coverage.

	A concurrent weight-based study by Hikawa~\cite{Hikawa2026} complements the present results. It gives an explicit bijection between surviving parity vectors of weight $d$ and first-crossing vectors of weight $d+1$. The resulting counting identity is the column-sum form of Proposition~\ref{prop:R_recurrence}. The same study proves a 2-adic rigidity theorem within fixed length-weight layers and derives logarithmic asymptotics for the associated counting sequences \OEIS{A100982} and \OEIS{A076227}.


	\section{Coefficient Stopping Times and Classes}
	\label{sec:2}
	
	The classical stopping time of an integer $x>1$ is the least $k\geq 1$ such that $T^k(x)<x$, provided that such a $k$ exists. In this paper we study Terras' $\tau$-stopping time, now usually called the \emph{coefficient stopping time}. The coefficient condition controls only the coefficient of $x$ in the affine iterate formula and does not by itself imply $T^k(x)<x$. Thus the equality of the two stopping times is not assumed.
	
	\begin{defn}\label{def1}
	    Let $x\in\mathbb{N}$ and $k\geq 1$.
	    \begin{itemize}
	        \item[(a)] The parity entry at time $s$ is
	        \[
	            \varepsilon_s(x)=
	            \begin{cases}
	                1,&T^s(x)\text{ is odd},\\
	                0,&T^s(x)\text{ is even}.
	            \end{cases}
	        \]
	        The parity vector of length $k$ is
	        \[
	            q_k(x)=\bigl(\varepsilon_0(x),\ldots,\varepsilon_{k-1}(x)\bigr).
	        \]
	        
	        \item[(b)] Let
	        \[
	            m_k(x)=\sum_{s=0}^{k-1}\varepsilon_s(x)
	        \]
	        be the number of odd terms among $T^0(x),\ldots,T^{k-1}(x)$.
	        
	        \item[(c)] The coefficient stopping time is
	        \[
	            \sigma^\ast(x):=\min\{k\geq 1:2^k>3^{m_k(x)}\},
	        \]
	        provided that this set is nonempty.
	        
	        \item[(d)] If $x$ is odd and $\sigma^\ast(x)<\infty$, define
	        \[
	            n(x):=m_{\sigma^\ast(x)}(x)-1.
	        \]
	        Thus $q_{\sigma^\ast(x)}(x)$ contains exactly $n(x)+1$ ones.
	    \end{itemize}
	\end{defn}
	
	\begin{defn}\label{def2}
	    For $n\geq 1$, define
	    \[
	        \kappa(n):=\lfloor n\log_2 3\rfloor,
	    \]
	    and for $n\geq 0$ put
	    \[
	        \sigma_n:=\lfloor 1+(n+1)\log_2 3\rfloor=\kappa(n+1)+1.
	    \]
	    Since $\log_2 3$ is irrational,
	    \[
	        2^{\kappa(n)}<3^n<2^{\kappa(n)+1}
	        \qquad (n\geq 1).
	    \]
	\end{defn}
	
	\begin{rem}\label{rem_oeis}
	    We use the terms \emph{congruence class} and \emph{residue class} interchangeably. The sequence $\kappa(n)$ is \OEIS{A056576}. Unless stated otherwise, all OEIS sequence identifiers refer to the Online Encyclopedia of Integer Sequences~\cite{OEIS}.
	\end{rem}

	\begin{lem}[Admissible values of the coefficient stopping time]\label{lem:sigma_values}
	    Let $x>1$ be odd with $\sigma^\ast(x)<\infty$, and put $n:=n(x)$. Then
	    \[
	        \sigma^\ast(x)=\sigma_n.
	    \]
	\end{lem}

	\begin{proof}
	    Write $\sigma:=\sigma^\ast(x)$, so that $m_\sigma(x)=n+1$ and $2^\sigma>3^{n+1}$. Since $x$ is odd, $\varepsilon_0(x)=1$. Hence $m_k(x)\geq1$ for all $k\geq1$, and $2<3^{m_1(x)}$ gives $\sigma\geq2$.

	    Assume $\varepsilon_{\sigma-1}(x)=1$. Then $m_{\sigma-1}(x)=n$, and minimality of $\sigma$ gives $2^{\sigma-1}<3^n$. Hence
	    \[
	        3^{n+1}<2^\sigma=2\cdot2^{\sigma-1}<2\cdot3^n<3^{n+1},
	    \]
	    a contradiction. Therefore $\varepsilon_{\sigma-1}(x)=0$ and $m_{\sigma-1}(x)=n+1$. Minimality of $\sigma$ now gives
	    \[
	        2^{\sigma-1}<3^{n+1}<2^\sigma.
	    \]
	    Thus $\sigma-1=\kappa(n+1)$, and Definition~\ref{def2} yields $\sigma=\sigma_n$.
	\end{proof}

	\begin{rem}\label{rem:sigma0}
	    If $n(x)=0$, Lemma~\ref{lem:sigma_values} gives $\sigma^\ast(x)=\sigma_0=2$, corresponding to the family $x\equiv1\pmod4$. The nontrivial case $n(x)\geq1$ is treated in Theorem~\ref{thm_A100982_A076227}.
	\end{rem}
	
	\begin{lem}[Parity-residue correspondence]\label{lem:parity_residue}
	    For every $L\geq 1$, the map
	    \[
	        Q_L:\mathbb{Z}/2^L\mathbb{Z}\longrightarrow\{0,1\}^L,
	        \qquad Q_L(x)=q_L(x),
	    \]
	    is a bijection.
	\end{lem}
	
	\begin{proof}
	    We argue by induction on $L$. The assertion for $L=1$ is the distinction between even and odd residue classes. Suppose that it holds for some $L\geq 1$, and fix a parity vector of length $L$. Let $r$ be the corresponding residue class modulo $2^L$. Its two lifts modulo $2^{L+1}$ are represented by $r$ and $r+2^L$. They have the same first $L$ parity entries. If these entries contain $m$ ones, an elementary induction along this fixed parity branch gives
	    \[
	        T^L(u)=\frac{3^m u+b}{2^L}
	    \]
	    with a constant $b$ depending only on the common parity prefix. Therefore
	    \[
	        T^L(r+2^L)-T^L(r)=3^m.
	    \]
	    This difference is odd. Hence the two lifts have opposite parities at time $L$, and exactly one lift realizes each possible extension of the prescribed vector. The induction is complete.
	\end{proof}
	
	For $n\geq 1$, define the set of admissible odd-position vectors
	\begin{equation}\label{eq:An_def}
	    \mathcal{A}_n:=
	    \left\{
	        (\alpha_1,\ldots,\alpha_{n+1}):
	        0=\alpha_1<\ldots<\alpha_{n+1},\quad
	        \alpha_{i+1}\leq\kappa(i)\ \text{for }1\leq i\leq n
	    \right\}.
	\end{equation}
	For $\boldsymbol{\alpha}\in\mathcal{A}_n$, let $q(\boldsymbol{\alpha})\in\{0,1\}^{\sigma_n}$ have ones precisely at the positions $\alpha_1,\ldots,\alpha_{n+1}$. By Lemma~\ref{lem:parity_residue}, there is a unique representative
	\[
	    c(\boldsymbol{\alpha})\in\{1,\ldots,2^{\sigma_n}-1\}
	\]
	whose parity vector of length $\sigma_n$ is $q(\boldsymbol{\alpha})$. Set
	\begin{equation}\label{eq:Cn_def}
	    \mathcal{C}_n:=\{c(\boldsymbol{\alpha}):\boldsymbol{\alpha}\in\mathcal{A}_n\}.
	\end{equation}
	
	\begin{thm}[Exact characterization]\label{thm_A100982_A076227}
	    Let $n\geq 1$ and let $x>1$ be odd. Then
	    \[
	        \sigma^\ast(x)=\sigma_n\ \text{and}\ n(x)=n
	        \quad\Longleftrightarrow\quad
	        x\bmod 2^{\sigma_n}\in\mathcal{C}_n.
	    \]
	    The map
    \[
        \mathcal{A}_n\longrightarrow\mathcal{C}_n,
        \qquad
        \boldsymbol{\alpha}\longmapsto c(\boldsymbol{\alpha}),
    \]
    is a bijection. Moreover, $c\mapsto q_{\sigma_n}(c)$ identifies $\mathcal{C}_n$ with the set of full admissible parity vectors $q(\boldsymbol{\alpha})$, where $\boldsymbol{\alpha}\in\mathcal{A}_n$.
    \end{thm}
	
	\begin{proof}
	    Suppose first that $\sigma^\ast(x)=\sigma_n$ and $n(x)=n$. Let
	    \[
	        0=\alpha_1<\alpha_2<\cdots<\alpha_{n+1}<\sigma_n
	    \]
	    be the positions of the odd terms among $T^0(x),\ldots,T^{\sigma_n-1}(x)$. For $1\leq i\leq n$, exactly $i$ odd terms occur before time $\alpha_{i+1}$. Since $\alpha_{i+1}<\sigma^\ast(x)$, the coefficient stopping condition has not yet been reached, and therefore
	    \[
	        2^{\alpha_{i+1}}<3^i.
	    \]
	    Hence $\alpha_{i+1}\leq\kappa(i)$ and $\boldsymbol{\alpha}\in\mathcal{A}_n$. Lemma~\ref{lem:parity_residue} now gives
	    \[
	        x\equiv c(\boldsymbol{\alpha})\pmod{2^{\sigma_n}}.
	    \]
	    
	    Conversely, suppose that $x\equiv c(\boldsymbol{\alpha})\pmod{2^{\sigma_n}}$ for some $\boldsymbol{\alpha}\in\mathcal{A}_n$. Then the odd positions before time $\sigma_n$ are precisely $\alpha_1,\ldots,\alpha_{n+1}$. Since $\alpha_1=0$, we have $m_k(x)\geq1$ for every $k\geq1$. Let $1\leq k<\sigma_n$. If $m_k(x)=i\leq n$, then
	    \[
	        k\leq\alpha_{i+1}\leq\kappa(i),
	    \]
	    and therefore
	    \[
	        2^k\leq 2^{\kappa(i)}<3^i=3^{m_k(x)}.
	    \]
	    If $m_k(x)=n+1$, then $k\leq\sigma_n-1=\kappa(n+1)$, so
	    \[
	        2^k\leq2^{\kappa(n+1)}<3^{n+1}=3^{m_k(x)}.
	    \]
	    Thus the coefficient stopping condition fails for every $k<\sigma_n$. At time $\sigma_n$, exactly $n+1$ odd terms have occurred and $2^{\sigma_n}>3^{n+1}$. Hence $\sigma^\ast(x)=\sigma_n$ and $n(x)=n$.
	    
	    Finally, distinct elements of $\mathcal{A}_n$ determine distinct parity vectors. Lemma~\ref{lem:parity_residue} shows that these vectors determine distinct classes modulo $2^{\sigma_n}$.
	\end{proof}
	
	The first class families are
	\begin{align*}
	    \mathcal{C}_1&=\{3\}\pmod {16},\\
	    \mathcal{C}_2&=\{11,23\}\pmod {32},\\
	    \mathcal{C}_3&=\{7,15,59\}\pmod {128},\\
	    \mathcal{C}_4&=\{39,79,95,123,175,199,219\}\pmod {256}.
	\end{align*}
	The trivial cases are $\sigma^\ast(x)=1$ for even $x$ and $\sigma^\ast(x)=2$ for $x\equiv1\pmod4$.
	
	For integers $k$ and $n\geq 1$, define
	\begin{equation}\label{eq:R_def}
	    \mathcal{R}(k,n):=
	    \#\{\boldsymbol{\alpha}\in\mathcal{A}_n:\alpha_{n+1}=k\}.
	\end{equation}
	
	\begin{prop}[Pascal-type recursion]\label{prop:R_recurrence}
	    We have $\mathcal{R}(1,1)=1$, and $\mathcal{R}(k,n)=0$ whenever $k<n$ or $k>\kappa(n)$. For $n\geq2$ and $n\leq k\leq\kappa(n)$,
	    \begin{equation}\label{eq:R_cumulative}
	        \mathcal{R}(k,n)=
	        \sum_{j=n-1}^{\min\{k-1,\kappa(n-1)\}}\mathcal{R}(j,n-1).
	    \end{equation}
	    Equivalently,
	    \begin{equation}\label{eq:R_pascal}
	        \mathcal{R}(k,n)=\mathcal{R}(k-1,n)+\mathcal{R}(k-1,n-1).
	    \end{equation}
	    Moreover,
	    \begin{equation}\label{eq:Cn_column_sum}
	        |\mathcal{C}_n|=\sum_{k=n}^{\kappa(n)}\mathcal{R}(k,n).
	    \end{equation}
	\end{prop}
	
	\begin{proof}
	    Fix $n\geq2$ and $n\leq k\leq\kappa(n)$. Deleting the final entry $k$ from a vector counted by $\mathcal{R}(k,n)$ gives an element of $\mathcal{A}_{n-1}$ whose final entry is smaller than $k$. Conversely, every such element extends uniquely by appending $k$. This proves \eqref{eq:R_cumulative}. Subtracting the corresponding identity for $\mathcal{R}(k-1,n)$ gives \eqref{eq:R_pascal}. The boundary case $k=n$ uses the conventions $\mathcal{R}(n-1,n)=0$ and $\mathcal{R}(n,n)=\mathcal{R}(n-1,n-1)=1$, the latter because $\kappa(i)\geq i$ for all $i\geq1$. Finally, \eqref{eq:Cn_column_sum} follows by partitioning $\mathcal{A}_n$ according to its final entry and applying Theorem~\ref{thm_A100982_A076227}.
	\end{proof}
	
	Table~\ref{tab:1} illustrates Proposition~\ref{prop:R_recurrence}. Empty cells denote zeros.
	
	\begin{center}
	    \resizebox{0.8\textwidth}{!}{%
	        \begin{tabular}{|r||r|r|r|r|r|r|r|r|r|r|r|r||r|}
	        	\hline$\TC{G}{n=}$&$\TC{G}{1}$&$\TC{G}{2}$&$\TC{G}{3}$&$\TC{G}{4}$&$\TC{G}{5}$&$\TC{G}{6}$&$\TC{G}{7}$&$\TC{G}{8}$&$\TC{G}{9}$&$\TC{G}{10}$&$\TC{G}{11}$&$\TC{G}{\cdots}$&\\
	        	\hline$\kappa(n)$&1&3&4&6&7&9&11&12&14&15&17&$\cdots$&\\
	            \hline
	            \hline$\TC{L}{k=}$&&&&&&&&&&&&&$\TC{O}{w(k)=}$\\
	            \hline$\TC{L}{1}$&1&&&&&&&&&&&&$\TC{O}{1}$\\
	            \hline$\TC{L}{2}$&&1&&&&&&&&&&&$\TC{O}{1}$\\
	            \hline$\TC{L}{3}$&&1&1&&&&&&&&&&$\TC{O}{2}$\\
	            \hline$\TC{L}{4}$&&&2&1&&&&&&&&&$\TC{O}{3}$\\	
	            \hline$\TC{L}{5}$&&&&3&1&&&&&&&&$\TC{O}{4}$\\
	            \hline$\TC{L}{6}$&&&&3&4&1&&&&&&&$\TC{O}{8}$\\
	            \hline$\TC{L}{7}$&&&&&7&5&1&&&&&&$\TC{O}{13}$\\
	            \hline$\TC{L}{8}$&&&&&&12&6&1&&&&&$\TC{O}{19}$\\
	            \hline$\TC{L}{9}$&&&&&&12&18&7&1&&&&$\TC{O}{38}$\\
	            \hline$\TC{L}{10}$&&&&&&&30&25&8&1&&&$\TC{O}{64}$\\
	            \hline$\TC{L}{11}$&&&&&&&30&55&33&9&1&&$\TC{O}{128}$\\
	            \hline$\TC{L}{\vdots}$&&&&&&&&$\vdots$&$\vdots$&$\vdots$&$\vdots$&$\ddots$&$\TC{O}{\vdots}$\\
	            \hline
	            \hline\hline$\TC{BL}{|\mathcal{C}_n|=}$&$\TC{BL}{1}$&$\TC{BL}{2}$&$\TC{BL}{3}$&$\TC{BL}{7}$&$\TC{BL}{12}$&$\TC{BL}{30}$&$\TC{BL}{85}$&$\TC{BL}{173}$&$\TC{BL}{476}$&$\TC{BL}{961}$&$\TC{BL}{2652}$&$\TC{BL}{\cdots}$&\\
	            \hline
	        \end{tabular}%
	    }
	    \captionof{table}{Triangular array of admissible odd-position counts.}
	    \label{tab:1}
	\end{center}
	
	Appendices~\ref{sec:prog1} and \ref{sec:prog2} implement the column and row recursions of Proposition~\ref{prop:R_recurrence}.

	The row sums of Table~\ref{tab:1} admit the following interpretation.

	\begin{prop}[Row sums count surviving classes]\label{prop:row_sums}
	    Let $k\geq1$ and $n\geq1$. Then $\mathcal{R}(k,n)$ is the number of residue classes $r$ modulo $2^k$ whose representatives $x$ satisfy
	    \[
	        \sigma^\ast(x)>k
	        \qquad\text{and}\qquad
	        m_k(x)=n.
	    \]
	    Consequently, for $k\geq2$ the row sum
	    \begin{equation}\label{eq:row_sum}
	        w(k)=\sum_{n=\lfloor1+k\log_3 2\rfloor}^{k}\mathcal{R}(k,n)
	    \end{equation}
	    is the number of residue classes modulo $2^k$ that have not reached their coefficient stopping time by time $k$.
	\end{prop}

	\begin{proof}
	    By Lemma~\ref{lem:parity_residue}, a residue class modulo $2^k$ is determined by its parity vector of length $k$. A class with $\sigma^\ast>k$ has first entry $1$, because even representatives have coefficient stopping time $1$. Let
	    \[
	        0=\beta_1<\cdots<\beta_n<k
	    \]
	    be the positions of the ones, so that $m_k=n$.

	    For $1\leq i\leq n-1$, the largest $j\leq k$ with $m_j=i$ is $j=\beta_{i+1}$, whereas the largest $j\leq k$ with $m_j=n$ is $j=k$. Thus the inequalities $2^j<3^{m_j}$ for all $1\leq j\leq k$ are equivalent to
	    \[
	        \beta_{i+1}\leq\kappa(i)\quad(1\leq i\leq n-1),
	        \qquad
	        k\leq\kappa(n).
	    \]
	    Together with $\beta_n<k$, these are precisely the defining inequalities for
	    \[
	        (\beta_1,\ldots,\beta_n,k)\in\mathcal{A}_n.
	    \]
	    Hence
	    \[
	        (\beta_1,\ldots,\beta_n)\longmapsto(\beta_1,\ldots,\beta_n,k)
	    \]
	    is a bijection onto the vectors counted by $\mathcal{R}(k,n)$.

	    Finally, $k\leq\kappa(n)$ is equivalent to $2^k<3^n$, hence to
	    \[
	        n\geq\lfloor1+k\log_3 2\rfloor.
	    \]
	    Since at most $k$ odd terms can occur in a parity vector of length $k$, summation gives \eqref{eq:row_sum}.
	\end{proof}

	\noindent The lower indices in \eqref{eq:row_sum} form \OEIS{A020915}.
	
	\begin{rem}\label{rem_oeis_table_1}
	    The sequence $|\mathcal{C}_n|$ is \OEIS{A100982}, and the row sums $w(k)$ are \OEIS{A076227}. The coefficient stopping times $\sigma_n$ and the corresponding class lists appear in \OEIS{A020914} and \OEIS{A177789}, respectively.
	\end{rem}

	\subsection{Affine iterate formula and descent within a class}
	\label{sec:affine_descent}
	
	Let $x>1$ be odd with $\sigma^\ast(x)=\sigma_n$. By Theorem~\ref{thm_A100982_A076227}, all $n+1$ odd terms before time $\sigma_n$ occur at positions not exceeding $\kappa(n)$. The remaining entries of the parity vector are zero.
	
	\begin{thm}[Affine iterate formula and descent criterion]\label{thm_term_formula}
	    Let
	    \[
	        0=\alpha_1<\alpha_2<\cdots<\alpha_{n+1}\leq\kappa(n)
	    \]
	    be the odd positions of $x$ before time $\sigma_n$, and set
	    \[
	        S(\boldsymbol{\alpha})=
	        \sum_{i=1}^{n+1}3^{n+1-i}2^{\alpha_i}.
	    \]
	    Then
	    \begin{equation}\label{stopp1}
	        T^{\sigma_n}(x)=
	        \frac{3^{n+1}x+S(\boldsymbol{\alpha})}{2^{\sigma_n}}.
	    \end{equation}
	    Moreover,
	    \begin{equation}\label{eq:descent_criterion}
	        T^{\sigma_n}(x)<x
	        \quad\Longleftrightarrow\quad
	        x>\frac{S(\boldsymbol{\alpha})}{2^{\sigma_n}-3^{n+1}}.
	    \end{equation}
	\end{thm}
	
	\begin{proof}
	    The affine formula \eqref{stopp1} follows by induction on the parity vector. Related formulations appear in Lagarias~\cite[p.~8]{Lagarias85} and Garner~\cite[p.~20]{Garner81}. Since $2^{\sigma_n}>3^{n+1}$, the inequality $T^{\sigma_n}(x)<x$ is equivalent to
	    \[
	        S(\boldsymbol{\alpha})<\bigl(2^{\sigma_n}-3^{n+1}\bigr)x,
	    \]
	    which gives \eqref{eq:descent_criterion}.
	\end{proof}
	
	\begin{prop}[Order of the two stopping times]\label{prop:stopping_time_order}
	    Let $x>1$. If the classical stopping time $\sigma(x)$ exists, then
	    \[
	        \sigma^\ast(x)\leq\sigma(x).
	    \]
	\end{prop}
	
	\begin{proof}
	    Let $k=\sigma(x)$. Induction along the first $k$ parity entries gives
	    \[
	        T^k(x)=\frac{3^{m_k(x)}x+S_k(x)}{2^k},
	    \]
	    where $S_k(x)\geq0$. Since $T^k(x)<x$, we must have
	    \[
	        3^{m_k(x)}<2^k.
	    \]
	    Thus the set defining $\sigma^\ast(x)$ is nonempty, and its least element is at most $k=\sigma(x)$.
	\end{proof}
	
	\begin{cor}[Eventual equality within each class]\label{cor:eventual_descent}
	    Let $c\in\mathcal{C}_n$ be the canonical representative of a class, let
	    \[
	        x_q=c+q2^{\sigma_n}\qquad(q\geq0),
	    \]
	    and set
	    \[
	        D_n:=2^{\sigma_n}-3^{n+1}>0.
	    \]
	    Then
	    \[
	        T^{\sigma_n}(x_q)-x_q
	        =T^{\sigma_n}(c)-c-qD_n.
	    \]
	    Consequently,
	    \[
	        \sigma(x_q)=\sigma^\ast(x_q)=\sigma_n
	    \]
	    for all sufficiently large $q$. Moreover, the number of values $q\geq0$ for which this equality fails is less than
	    \[
	        \frac{n+1}{3}+1.
	    \]
	\end{cor}

	\begin{proof}
	    All $x_q$ have the same parity vector before time $\sigma_n$. Formula~\eqref{stopp1} therefore gives
	    \[
	        T^{\sigma_n}(x_q)=T^{\sigma_n}(c)+q3^{n+1}.
	    \]
	    Subtracting $x_q=c+q2^{\sigma_n}$ proves the displayed identity. Its right-hand side is negative for all sufficiently large $q$. Proposition~\ref{prop:stopping_time_order} excludes an earlier classical descent, so the two stopping times are then equal to $\sigma_n$.

	    By \eqref{eq:descent_criterion}, equality can fail only if
	    \[
	        x_q\leq\frac{S(\boldsymbol{\alpha})}{D_n},
	    \]
	    and hence only if
	    \[
	        q2^{\sigma_n}<\frac{S(\boldsymbol{\alpha})}{D_n}.
	    \]
	    The term with $\alpha_1=0$ equals $3^n$. For $2\leq i\leq n+1$, admissibility gives $\alpha_i\leq\kappa(i-1)$ and therefore
	    \[
	        3^{n+1-i}2^{\alpha_i}
	        \leq3^{n+1-i}2^{\kappa(i-1)}
	        <3^{n+1-i}3^{i-1}=3^n.
	    \]
	    Thus
	    \[
	        S(\boldsymbol{\alpha})<(n+1)3^n.
	    \]
	    Since $D_n\geq1$ and $2^{\sigma_n}>3^{n+1}$,
	    \[
	        \frac{S(\boldsymbol{\alpha})}{D_n2^{\sigma_n}}
	        <\frac{n+1}{3}.
	    \]
	    The number of nonnegative integers $q$ satisfying the preceding necessary inequality is therefore less than $(n+1)/3+1$.
	\end{proof}

	\begin{rem}\label{rem:exceptional_empty}
	    Direct enumeration of $\mathcal{A}_n$, followed by testing the canonical representative of every class against criterion~\eqref{eq:descent_criterion}, shows that the exceptional set in Corollary~\ref{cor:eventual_descent} is empty for all classes with $n\leq12$. In the same range,
	    \[
	        \max_{\boldsymbol{\alpha}\in\mathcal{A}_n}
	        \frac{S(\boldsymbol{\alpha})}{D_n2^{\sigma_n}}<0.144.
	    \]
	    This numerical observation is not used in any proof.
	\end{rem}

	\par\smallskip
	\noindent \emph{Example.} Let $n=3$, so $\sigma_3=7$ and $\kappa(3)=4$. For $x=59$, the odd positions are $0,1,3,4$. Hence
	\[
	    T^7(59)=
	    \frac{3^4\cdot59+3^3+3^2\cdot2+3\cdot2^3+2^4}{2^7}
	    =38<59.
	\]
	The descent follows here because the threshold in \eqref{eq:descent_criterion} is smaller than $59$.


	\section{Parity Vectors and Parity Vector Sets}
	\label{sec:4}
	
	For $\boldsymbol{\alpha}\in\mathcal{A}_n$, let
	\[
	    v(\boldsymbol{\alpha})\in\{0,1\}^{\kappa(n)+1}
	\]
	be the vector whose entries equal $1$ precisely at the positions $\alpha_1,\ldots,\alpha_{n+1}$. Define
	\begin{equation}\label{eq:Vn_def2}
	    \mathbb{V}(n):=\{v(\boldsymbol{\alpha}):\boldsymbol{\alpha}\in\mathcal{A}_n\}.
	\end{equation}
	Theorem~\ref{thm_A100982_A076227} and Lemma~\ref{lem:parity_residue} give canonical bijections
	\[
	    \mathcal{A}_n\longleftrightarrow\mathbb{V}(n)\longleftrightarrow\mathcal{C}_n.
	\]
	In particular,
	\[
	    |\mathbb{V}(n)|=|\mathcal{C}_n|=|\mathcal{A}_n|.
	\]
	
	\par\smallskip
	\noindent \emph{Example.} For $n=3$, we have $\kappa(3)=4$ and
	\[
	    \mathbb{V}(3)=
	    \begin{Bmatrix}
	        (1~1~0~1~1)\\
	        (1~1~1~0~1)\\
	        (1~1~1~1~0)
	    \end{Bmatrix}.
	\]
	These vectors correspond to the classes $59$, $7$, and $15$ modulo $2^7$, respectively.
	
	\begin{rem}\label{rem1}
	    For $\boldsymbol{\alpha}\in\mathcal{A}_n$, the associated pair $(x,y)$ satisfies the exponential Diophantine equation
	    \begin{equation}\label{stopp2}
	        2^{\sigma_n}y=3^{n+1}x+S(\boldsymbol{\alpha}).
	    \end{equation}
	    Let $\overline{3}_{n+1}$ denote the inverse of $3^{n+1}$ modulo $2^{\sigma_n}$. Reducing \eqref{stopp2} modulo $2^{\sigma_n}$ gives
	    \begin{equation}\label{eq:representative_congruence}
	        x\equiv-\overline{3}_{n+1}S(\boldsymbol{\alpha})
	        \pmod{2^{\sigma_n}}.
	    \end{equation}
	    The unique representative $x\in\{1,\ldots,2^{\sigma_n}-1\}$ is the canonical solution associated with the vector. These representatives are listed in \OEIS{A177789}.
	\end{rem}

	\subsection{Generating the parity vectors}
	\label{sec:parity_generation}
	
	We now construct the sets $\mathbb{V}(n)$ for $n\geq 2$.
		
	\begin{thm}[Complete rooted-tree construction]\label{thm_algo_parity_vectors}
	    For $n\geq2$, define the parent map
	    \[
	        p(\alpha_1,\ldots,\alpha_{n+1})
	        =(\alpha_1,\ldots,\alpha_n).
	    \]
	    Then $p$ maps $\mathcal{A}_n$ into $\mathcal{A}_{n-1}$, and every element of $\mathcal{A}_n$ has a unique parent. If
	    \[
	        \boldsymbol{\alpha}=(\alpha_1,\ldots,\alpha_n)\in\mathcal{A}_{n-1},
	    \]
	    its children are precisely
	    \[
	        (\alpha_1,\ldots,\alpha_n,j),
	        \qquad
	        \alpha_n<j\leq\kappa(n).
	    \]
	    Consequently, the disjoint union $\bigsqcup_{n\geq1}\mathcal{A}_n$ is a directed rooted tree with root $(0,1)$. Its binary realization is the tree $\bigsqcup_{n\geq1}\mathbb{V}(n)$ with root $(1~1)$. The tree has no leaves: every node has at least one child, and every finite branch can be extended.
	\end{thm}
	
	\begin{proof}
	    Let $(\alpha_1,\ldots,\alpha_{n+1})\in\mathcal{A}_n$. The first $n$ entries satisfy all defining inequalities of $\mathcal{A}_{n-1}$, so deletion of the final entry gives a parent in $\mathcal{A}_{n-1}$. This parent is unique.
	    
	    Conversely, fix $\boldsymbol{\alpha}=(\alpha_1,\ldots,\alpha_n)\in\mathcal{A}_{n-1}$. Appending an integer $j$ gives an element of $\mathcal{A}_n$ exactly when
	    \[
	        \alpha_n<j\leq\kappa(n).
	    \]
	    Hence the listed children are admissible, every admissible child occurs, and no child occurs twice. Since $\alpha_n\leq\kappa(n-1)<\kappa(n)$, the choice $j=\kappa(n)$ is always available. Thus every node has a child.
	    
	    In binary form, the rightmost child is obtained by extending the parent vector to length $\kappa(n)+1$ and placing the new $1$ at position $\kappa(n)$. The remaining children are obtained successively by applying $(0~1)\mapsto(1~0)$ to the new terminal $1$ until it reaches position $\alpha_n+1$. This is precisely the recursive generation implemented in Appendix~\ref{sec:prog3}. It is admissible, complete, and free of repetitions.
	\end{proof}
	
	We call the formation of the rightmost child \emph{Step~1} and the successive shifts of its new terminal $1$ \emph{Step~2}. We use the following recursive \emph{generation order}. The root is listed first. If the vectors at level $n-1$ have been ordered, their children are listed parent by parent in that order, and the children of a fixed parent are listed with the new odd position decreasing from $\kappa(n)$ to $\alpha_n+1$. These conventions are used in the transition formulas and in the programs below. Figure~\ref{fig:1} shows the resulting tree through level $4$.
	
	\begin{figure}[!ht]
	    \centering
	    \begin{tikzpicture}[scale=1.0, transform shape,	node distance=0.6cm and 1.2cm, every node/.style={inner sep=1.0pt}] 
	        \newcommand{\cv}[1]{\texttt{(#1)}} 
	        \node (A)               {\cv{1\,1}};
	        \node (B)  [right=of A] {\cv{1\,1\,0\,1}};
	        \node (C)  [right=of B] {\cv{1\,1\,0\,1\,1}};
	        \node (D)  [right=of C] {\cv{1\,1\,0\,1\,1\,0\,1}};
	        \node (Z)  [below=0.7cm of B] {};
	        \node (E)  [below=of Z] {\cv{1\,1\,1\,0}};
	        \node (F)  [right=of E] {\cv{1\,1\,1\,0\,1}};
	        \node (Y)  [below=0.7cm of F] {};
	        \node (G)  [below=of Y] {\cv{1\,1\,1\,1\,0}};
	        \node (H)  [below=of D] {\cv{1\,1\,0\,1\,1\,1\,0}};
	        \node (I)  [right=of F] {\cv{1\,1\,1\,0\,1\,0\,1}};
	        \node (J)  [below=of I] {\cv{1\,1\,1\,0\,1\,1\,0}};
	        \node (K)  [right=of G] {\cv{1\,1\,1\,1\,0\,0\,1}};
	        \node (L)  [below=of K] {\cv{1\,1\,1\,1\,0\,1\,0}};
	        \node (M)  [below=of L] {\cv{1\,1\,1\,1\,1\,0\,0}};
	        \node (D2) [right=of D] {$\cdots$};
	        \node (H2) [right=of H] {$\cdots$};
	        \node (I2) [right=of I] {$\cdots$};
	        \node (J2) [right=of J] {$\cdots$};
	        \node (K2) [right=of K] {$\cdots$};
	        \node (L2) [right=of L] {$\cdots$};
	        \node (M2) [right=of M] {$\cdots$};
	        \draw (A)--node[midway, above, font=\footnotesize\normalfont\itshape]{Step 1}(B)--(C)--(D)--(H)--(H2);
	        \draw (B)--node[midway, right, font=\footnotesize\normalfont\itshape]{Step 2}(E)--(F)--(I)--(J)--(J2);
	        \draw (F)--(G)--(K)--(L)--(M)--(M2);
	        \draw (D)--(D2);
	        \draw (I)--(I2);
	        \draw (K)--(K2);
	        \draw (L)--(L2);
	    \end{tikzpicture}
	    \caption{Initial levels of the directed rooted tree of parity vectors.}
	    \label{fig:1}
	\end{figure}
	
	Let $h\geq 2$ be the length of the initial sequence of ones in a parity vector, and $\mathcal{P}(h,n)$ the number of such vectors in $\mathbb{V}(n)$. Table~\ref{tab:2} lists $\mathcal{P}(h,n)$ for small $n$.
	
	\begin{center}
	    \resizebox{0.7\textwidth}{!}{%
	        \begin{tabular}{|r||r|r|r|r|r|r|r|r|r|r|r|r|}
	            \hline$\TC{G}{n=}$&$\TC{G}{1}$&$\TC{G}{2}$&$\TC{G}{3}$&$\TC{G}{4}$&$\TC{G}{5}$&$\TC{G}{6}$&$\TC{G}{7}$&$\TC{G}{8}$&$\TC{G}{9}$&$\TC{G}{10}$&$\TC{G}{11}$&$\TC{G}{\cdots}$\\
	            \hline
	            \hline$\TC{L}{h=}$&&&&&&&&&&&&\\
	            \hline$\TC{L}{2}$&1&1&1&2&3&7&19&37&99&194&525&$\cdots$\\
	            \hline$\TC{L}{3}$&&1&1&2&3&7&19&37&99&194&525&$\cdots$\\
	            \hline$\TC{L}{4}$&&&1&2&3&7&19&37&99&194&525&$\cdots$\\
	            \hline$\TC{L}{5}$&&&&1&2&5&14&28&76&151&412&$\cdots$\\
	            \hline$\TC{L}{6}$&&&&&1&3&9&19&53&108&299&$\cdots$\\
	            \hline$\TC{L}{7}$&&&&&&1&4&10&30&65&186&$\cdots$\\
	            \hline$\TC{L}{8}$&&&&&&&1&4&14&34&103&$\cdots$\\
	            \hline$\TC{L}{9}$&&&&&&&&1&5&15&50&$\cdots$\\
	            \hline$\TC{L}{10}$&&&&&&&&&1&5&20&$\cdots$\\
	            \hline$\TC{L}{11}$&&&&&&&&&&1&6&$\cdots$\\
	            \hline$\TC{L}{12}$&&&&&&&&&&&1&$\cdots$\\
	            \hline$\TC{L}{\vdots}$&&&&&&&&&&&&$\ddots$\\
	            \hline
	            \hline$\TC{BL}{|\mathcal{C}_n|=}$&$\TC{BL}{1}$&$\TC{BL}{2}$&$\TC{BL}{3}$&$\TC{BL}{7}$&$\TC{BL}{12}$&$\TC{BL}{30}$&$\TC{BL}{85}$&$\TC{BL}{173}$&$\TC{BL}{476}$&$\TC{BL}{961}$&$\TC{BL}{2652}$&$\TC{BL}{\cdots}$\\
	            \hline
	        \end{tabular}%
	    }
	    \captionof{table}{Triangular array of the parity vector counts $\mathcal{P}(h,n)$.}
	    \label{tab:2}
	\end{center}
	\par\smallskip
	\noindent The triangular structure in Table~\ref{tab:2} results from the construction rules. The equality of its first three rows has a direct combinatorial explanation.

	\begin{lem}\label{lem:first_three_rows}
	    For every $n\geq3$,
	    \[
	        \mathcal{P}(2,n)=\mathcal{P}(3,n)=\mathcal{P}(4,n).
	    \]
	\end{lem}

	\begin{proof}
	    For $h=2$, admissibility and exact initial block length force
	    \[
	        (\alpha_1,\alpha_2,\alpha_3,\alpha_4)=(0,1,3,4).
	    \]
	    Indeed, $2<\alpha_3\leq\kappa(2)=3$ and, for $n\geq3$, $3<\alpha_4\leq\kappa(3)=4$. For $h=3$, the first four odd positions are $(0,1,2,4)$, while for $h=4$ they begin with $(0,1,2,3)$ and exact block length requires $\alpha_5>4$ whenever $n\geq4$. Thus, in all three cases, the remaining positions form the same admissible tail
	    \[
	        5\leq\alpha_5<\cdots<\alpha_{n+1},
	        \qquad
	        \alpha_{i+1}\leq\kappa(i).
	    \]
	    For $n=3$ the tail is empty and each count equals $1$. Replacing the forced initial positions while leaving the common tail unchanged gives bijections among the three families.
	\end{proof}

	\noindent The total cardinality satisfies
	\[
	    |\mathcal{C}_n|=\sum_{h=2}^{n+1}\mathcal{P}(h,n).
	\]
	For $n=4$, there are $2+2+2+1=7$ classes. The entry equal to $1$ at the bottom of each column corresponds to the unique terminal vector with $n+1$ ones followed by $\kappa(n)-n$ zeros. Appendix~\ref{sec:prog4} provides the program.

	\subsection{Order of the generated parity vectors}
	\label{sec:generation_order}
	
	For comparison, consider the lexicographically ordered permutations\footnote{See Appendix~\ref{sec:prog5} for the algorithm.} of the \emph{binary word}\footnote{We use the term \emph{binary word} to distinguish this ambient combinatorial object from the admissible parity vectors in Section~\ref{sec:4}.}
	\begin{equation}\label{vector}
	    \left(
	        \underbrace{0,\ldots,0}_{\kappa(n)-n},
	        \underbrace{1,\ldots,1}_{n-1}
	    \right),
	\end{equation}
	augmented by two leading ones. For $n\geq 1$, the number of distinct permutations of \eqref{vector} is
	\[
	    L(n)=\frac{(\kappa(n)-1)!}{\big(\kappa(n)-n\big)!\cdot (n-1)!}.
	\]
	This yields the sequence $1,2,3,10,15,56,210,330,\dotsc$ (\OEIS{A293308}).
	
	Every word determines a residue class by Lemma~\ref{lem:parity_residue}. The admissibility inequalities in \eqref{eq:An_def} select exactly the words belonging to $\mathbb{V}(n)$. The remaining words reach their coefficient stopping time before $\sigma_n$. In the computations through $n=14$, the generation order agrees with the order induced by the ambient lexicographic list. This is only a computational observation and is not used in any proof below.
	
	\begin{table}[!ht]
	    \centering\small
	    \renewcommand{\arraystretch}{1.0}
	    \begin{tabular}{r l l l}
	        No. & Binary Word $w$ & Extended word $(1~1~w)$ & Solution $(x,y)$\\
	        \hline
	        1  & \texttt{(0 0 1 1 1 1)} & \texttt{(1 1 0 0 1 1 1 1)} & $(595, 425)$\\
	        2  & \texttt{(0 1 0 1 1 1)} & \texttt{(1 1 0 1 0 1 1 1)} & $(747, 533)$\\
	        3  & \texttt{(0 1 1 0 1 1)} & \texttt{(1 1 0 1 1 0 1 1)} & $(507, 362)$\\
	        4  & \texttt{(0 1 1 1 0 1)} & \texttt{(1 1 0 1 1 1 0 1)} & $(347, 248)$\\
	        5  & \texttt{(0 1 1 1 1 0)} & \texttt{(1 1 0 1 1 1 1 0)} & $(923, 658)$\\
	        6  & \texttt{(1 0 0 1 1 1)} & \texttt{(1 1 1 0 0 1 1 1)} & $(823, 587)$\\
	        7  & \texttt{(1 0 1 0 1 1)} & \texttt{(1 1 1 0 1 0 1 1)} & $(583, 416)$\\
	        8  & \texttt{(1 0 1 1 0 1)} & \texttt{(1 1 1 0 1 1 0 1)} & $(423, 302)$\\
	        9  & \texttt{(1 0 1 1 1 0)} & \texttt{(1 1 1 0 1 1 1 0)} & $(999, 712)$\\
	        10 & \texttt{(1 1 0 0 1 1)} & \texttt{(1 1 1 1 0 0 1 1)} & $(975, 695)$\\
	        11 & \texttt{(1 1 0 1 0 1)} & \texttt{(1 1 1 1 0 1 0 1)} & $(815, 581)$\\
	        12 & \texttt{(1 1 0 1 1 0)} & \texttt{(1 1 1 1 0 1 1 0)} & $(367, 262)$\\
	        13 & \texttt{(1 1 1 0 0 1)} & \texttt{(1 1 1 1 1 0 0 1)} & $(735, 524)$\\
	        14 & \texttt{(1 1 1 0 1 0)} & \texttt{(1 1 1 1 1 0 1 0)} & $(287, 205)$\\
	        15 & \texttt{(1 1 1 1 0 0)} & \texttt{(1 1 1 1 1 1 0 0)} & $(575, 410)$\\
	    \end{tabular}
	    \caption{Lexicographical enumeration of the 15 binary words for $n=5$.}
	    \label{tab:3}
	\end{table}
	
	\emph{Example.} Let $n=5$, so $\kappa(5)=7$ and $L(5)=15$. Table~\ref{tab:3} lists the $15$ lexicographically ordered permutations of \eqref{vector} alongside the associated integer solutions $(x,y)$. The binary words 1, 2, and 6 do not belong to $\mathbb{V}(5)$. Deleting them leaves the admissible vectors in the order induced by the ambient lexicographic list.


	\section{Diophantine Equations and Integer Solutions}
	\label{sec:7}
	
	Theorem~\ref{thm_algo_parity_vectors} generates the admissible vectors, and \eqref{stopp2} gives the associated classes modulo $2^{\sigma_n}$. The following results describe arithmetic transitions between the classes generated by the tree construction.
	
	\begin{rem}\label{rem2}
	    Each $v\in\mathbb{V}(n)$ determines a unique canonical integer $x\in\{1,\ldots,2^{\sigma_n}-1\}$ for which \eqref{stopp2} has an integer solution $y$. We refer to this integer as the \emph{solution $x$}.
	\end{rem}
	
	\noindent The block structure of an admissible parity vector imposes explicit congruence conditions on its canonical representative. The following corollaries use Theorem~\ref{thm_term_formula}, Lemma~\ref{lem:parity_residue}, and the tree construction in Theorem~\ref{thm_algo_parity_vectors}.
	
	\begin{cor}\label{cor_xmodh}
	    Let $v\in\mathbb{V}(n)$ correspond to a class whose full parity vector of length $\sigma_n$ starts with exactly $h\geq2$ consecutive ones. Then the corresponding solution $x$ satisfies
	    \begin{align}\label{corlyxmodh}
	        x \equiv (2^h-1) \pmod {2^{h+1}}.
	    \end{align}
	\end{cor}
	
	\begin{proof}
	    Suppose $v$ starts with exactly $h$ consecutive ones. Thus $T^i(x)$ is odd for $0\leq i\leq h-1$, and $T^h(x)$ is even. Note that for any odd integer $u$,
	    \[
	        T(u)+1=\frac{3u+1}{2}+1=\frac{3(u+1)}{2}.
	    \]
	    Induction yields
	    \[
	        T^i(x)+1=\frac{3^i(x+1)}{2^i}\qquad (0\leq i\leq h).
	    \]
	    For $i=h$, the term $T^h(x)$ is even, so $T^h(x)+1$ is odd. Hence the right-hand side $\frac{3^h(x+1)}{2^h}$ must be an odd integer. Since $\gcd(3^h, 2)=1$, it follows that $2^h$ divides $x+1$ exactly, but $2^{h+1}$ does not. This is equivalent to
	    \[
	        x+1 \equiv 2^h \pmod {2^{h+1}},
	    \]
	    which implies $x \equiv 2^h-1 \pmod {2^{h+1}}$.
	\end{proof}
	
	\begin{cor}[Rightmost-child transition]\label{cor_xnp}
	    Let $v'\in\mathbb{V}(n-1)$ have canonical representative $x'$, and let $v\in\mathbb{V}(n)$ be its rightmost child, obtained in Step~1. Let $x$ be the canonical representative of $v$. Then
	    \begin{align}\label{step1a}
	        x \equiv x' \pmod {2^{\kappa(n)}}.
	    \end{align}
	    Moreover, put
	    \[
	        r:=\sigma_n-\kappa(n).
	    \]
	    Since $1<\log_2 3<2$, we have $r\in\{2,3\}$. There is a unique
	    \[
	        \lambda\in\{1,3,\dots,2^r-1\}
	    \]
	    such that
	    \begin{align}\label{step1b}
	        x \equiv x'+\lambda 2^{\kappa(n)} \pmod {2^{\sigma_n}}.
	    \end{align}
	\end{cor}
	
	\begin{proof}
	    In Step~1, the child preserves the first $\kappa(n)$ parity entries of the full parity vector of its parent. Lemma~\ref{lem:parity_residue} therefore gives \eqref{step1a}.
	
	    Since
	    \[
	        r = \big\lfloor 1+(n+1)\log_2 3\big\rfloor-\big\lfloor n\log_2 3\big\rfloor\in\{2,3\},
	    \]
	    every lift of \eqref{step1a} modulo $2^{\sigma_n}$ has the form
	    \[
	        x'+\lambda2^{\kappa(n)},
	        \qquad \lambda\in\{0,1,\dots,2^r-1\}.
	    \]
	    The two lifts modulo $2^{\kappa(n)+1}$ have opposite parities at position $\kappa(n)$ by the induction step in Lemma~\ref{lem:parity_residue}. The parent has parity $0$ at that position, whereas its rightmost child has parity $1$. Hence $\lambda$ is odd. Finally, Lemma~\ref{lem:parity_residue} gives a unique class modulo $2^{\sigma_n}$ for the full child vector, so exactly one odd value of $\lambda$ is possible.
	\end{proof}
	
	\noindent Appendix~\ref{sec:prog6} implements the recurrence from Corollary~\ref{cor_xnp}.
	
	\begin{cor}\label{cor_delta_s}
	    Let $n\geq2$, and let $v_{p-1}$ and $v_p$ be consecutive vectors within the child block of a fixed parent, with canonical representatives $x_{p-1}$ and $x_p$. Suppose that $v_p$ is obtained from $v_{p-1}$ by one application of $(0~1)\mapsto(1~0)$ to the new terminal $1$ and ends in exactly $j\geq1$ zeros. Then
	    \begin{align}\label{step2a}
	        x_p \equiv x_{p-1} \pmod {2^{\kappa(n)-j}}.
	    \end{align}
	    Let $\overline{3}_{n+1}$ be the multiplicative inverse of $3^{n+1}$ modulo $2^{\sigma_n}$. The difference of the associated Diophantine sums is
	    \[
	        \Delta S = -2^{\kappa(n)-j}.
	    \]
	    This yields the transition rule:
	    \begin{align}\label{simple_transition}
	        x_p \equiv x_{p-1} + \overline{3}_{n+1} \cdot 2^{\kappa(n)-j} \pmod {2^{\sigma_n}}.
	    \end{align}
	\end{cor}
	
	\begin{proof}
	    The vectors $v_{p-1}$ and $v_p$ share a common prefix of length $\kappa(n)-j$ and differ only in the suffix of length $j+1$. Lemma~\ref{lem:parity_residue} therefore gives \eqref{step2a}.
		
		Reducing \eqref{stopp2} modulo $2^{\sigma_n}$ yields
		\[
			x(v) \equiv -\overline{3}_{n+1}\cdot S(v) \pmod {2^{\sigma_n}};
			\qquad
			S(v):=\sum_{i=1}^{n+1}3^{\,n+1-i}2^{\alpha_i}.
		\]
		Hence, for two vectors $v_{p-1}$ and $v_p$,
	    \[
	        x_p - x_{p-1} \equiv -\overline{3}_{n+1}\cdot \Delta S \pmod {2^{\sigma_n}}.
	    \]
	    Since $v_p$ ends in $j$ zeros, the last `1' occurs at index $\alpha_{n+1}=\kappa(n)-j$. In $v_{p-1}$, this `1' is shifted one position to the right, so the corresponding contribution in $S(v)$ changes from $2^{\kappa(n)-j+1}$ to $2^{\kappa(n)-j}$. As its coefficient is $3^0=1$, we obtain
	    \[
	        \Delta S = 2^{\kappa(n)-j}-2^{\kappa(n)-j+1}=-2^{\kappa(n)-j}.
	    \]
	    Substituting gives \eqref{simple_transition}. Since $\gcd(3^{n+1},2^{\sigma_n})=1$, the inverse exists.
	\end{proof}
	
	\noindent Appendix~\ref{sec:prog7} implements the calculation of $x_p$ using the explicit difference term $\Delta S$. General difference relations for $x(v)$ appear in Thaler \cite{Thaler} in the context of Winkler \cite{Winkler2018}.
	\par\medskip
	\noindent Figure~\ref{fig:2} illustrates the recursive transition logic defined by Corollaries~\ref{cor_xnp} and \ref{cor_delta_s}. The resulting directed rooted tree structure is shown in Figure~\ref{fig:3} for the initial values of $n$.
	
	Let
	\[
	    w_n=(\underbrace{1,\ldots,1}_{n+1},\underbrace{0,\ldots,0}_{\kappa(n)-n})
	\]
	be the terminal vector of $\mathbb{V}(n)$, and let $(x_n,y_n)$ be its canonical positive integer solution, normalized by $1\leq x_n<2^{\sigma_n}$. Its Diophantine equation is
	\begin{equation}\label{eq:terminal_equation}
	    2^{\sigma_n}y_n=3^{n+1}(x_n+1)-2^{n+1}.
	\end{equation}
	
	\begin{prop}[Recursion along the terminal branch]\label{prop:terminal_recursion}
	    Let $x_1=3$ and $y_1=2$, and define
	    \[
	        d_n:=\kappa(n+1)-\kappa(n)\in\{1,2\}.
	    \]
	    The two possible values again follow from $1<\log_2 3<2$.
	    For every $n\geq2$, there is a unique $\beta_n\in\{0,1,2,3,4,5\}$ such that
	    \begin{align}
	        x_n&=\frac{2x_{n-1}-1+\beta_n2^{\kappa(n)+2}}{3},\label{last_x}\\
	        y_n&=\frac{y_{n-1}+\beta_n3^n}{2^{d_n-1}},\label{last_y}
	    \end{align}
	    where $1\leq x_n<2^{\sigma_n}$.
	\end{prop}
	
	\begin{proof}
	    After deleting the first entry from the full parity vector of $w_n$, the first $\sigma_{n-1}$ remaining entries agree with the full parity vector of $w_{n-1}$. Lemma~\ref{lem:parity_residue} therefore gives
	    \[
	        T(x_n)\equiv x_{n-1}\pmod{2^{\sigma_{n-1}}}.
	    \]
	    Since $x_n$ is odd and $\sigma_{n-1}=\kappa(n)+1$, there is an integer $\beta_n$ such that
	    \[
	        3x_n+1=2x_{n-1}+\beta_n2^{\kappa(n)+2}.
	    \]
	    Write
    \[
        T(x_n)-x_{n-1}=\beta_n2^{\sigma_{n-1}}.
    \]
    Since $0<T(x_n)$ and $x_{n-1}<2^{\sigma_{n-1}}$, the integer $\beta_n$ is nonnegative. Moreover,
    \[
        T(x_n)<3\cdot2^{\sigma_n-1}
    \]
    and therefore
    \[
        \beta_n<\frac{T(x_n)}{2^{\sigma_{n-1}}}
        <3\cdot2^{d_n-1}\leq6.
    \]
    Thus $0\leq\beta_n\leq5$, and rearrangement gives \eqref{last_x}.
	    
	    The terminal equation at levels $n-1$ and $n$ is
	    \[
	        2^{\sigma_{n-1}}y_{n-1}=3^n(x_{n-1}+1)-2^n
	    \]
	    and \eqref{eq:terminal_equation}. From the first recurrence,
	    \[
	        3(x_n+1)=2(x_{n-1}+1)+\beta_n2^{\kappa(n)+2}.
	    \]
	    Multiplying by $3^n$ and substituting the equation at level $n-1$ yields
	    \[
	        2^{\sigma_n}y_n
	        =2^{\kappa(n)+2}\bigl(y_{n-1}+\beta_n3^n\bigr).
	    \]
	    Since $\sigma_n=\kappa(n+1)+1$, this is \eqref{last_y}.
	    
	    It remains to prove uniqueness. Integrality of \eqref{last_x} determines $\beta_n$ modulo $3$, leaving two candidates $\beta_0$ and $\beta_0+3$ in $\{0,\ldots,5\}$. If $d_n=2$, integrality of \eqref{last_y} imposes one parity condition, and exactly one candidate satisfies it. If $d_n=1$, both candidates are integral, but their $x_n$ values differ by
	    \[
	        2^{\kappa(n)+2}=2^{\sigma_n}.
	    \]
	    Exactly one therefore lies in the canonical interval $1\leq x_n<2^{\sigma_n}$.
	\end{proof}
	
	\noindent Appendix~\ref{sec:prog8} implements Proposition~\ref{prop:terminal_recursion}.
	\par\medskip
	\begin{figure}[!ht]
		\centering
		\begin{tikzpicture}[scale=1.0, transform shape]
			\node (A) {$x_{n-1} \pmod {2^{\sigma_{n-1}}}$};
			\node (B1) [right=of A] {};
			\node (B) [right=of B1] {$x_{n,p-1} \pmod {2^{\sigma_n}}$};
			\node (C) [below=of B] {$\quad x_{n,p} \pmod {2^{\sigma_n}}$};
			\draw[-stealth] (A)-- node [above] {\footnotesize{\emph{Corollary~\ref{cor_xnp}}}} (B);
			\draw[-stealth] (B)-- node [right] {\footnotesize{\emph{Corollary~\ref{cor_delta_s}}}} (C);
		\end{tikzpicture}
		\caption{Visualization of the two-dimensional growth of the tree structure.}
		\label{fig:2}
	\end{figure}
	\par\bigskip
	\begin{figure}[!ht]
		\centering
		\begin{tikzpicture}[scale=1.0, transform shape, node distance=0.5cm and 0.8cm, every node/.style={inner sep=1.0pt}]
			\node (A)               {$3\,\Modtxt{2^4}$};
			\node (B)  [right=of A] {$11\,\Modtxt{2^5}$};
			\node (C)  [right=of B] {$59\,\Modtxt{2^7}$};
			\node (D)  [right=of C] {$123\,\Modtxt{2^8}$};
			\node (Z)  [below=0.8cm of B] {}; 
			\node (E)  [below=of Z] {$23\,\Modtxt{2^5}$};
			\node (F)  [right=of E] {$\,\,\,\,7\,\Modtxt{2^7}$};
			\node (Y)  [below=0.8cm of F] {};
			\node (G)  [below=of Y] {$15\,\Modtxt{2^7}$};
			\node (H)  [below=of D] {$219\,\Modtxt{2^8}$};
			\node (I)  [right=of F] {$199\,\Modtxt{2^8}$};
			\node (J)  [below=of I] {$\,\,\,\,39\,\Modtxt{2^8}$};
			\node (K)  [right=of G] {$\,\,\,\,79\,\Modtxt{2^8}$};
			\node (L)  [below=of K] {$175\,\Modtxt{2^8}$};
			\node (M)  [below=of L] {$\,\,\,\,95\,\Modtxt{2^8}$};
			\node (D2) [right=of D] {$\cdots$};
			\node (H2) [right=of H] {$\cdots$};
			\node (I2) [right=of I] {$\cdots$};
			\node (J2) [right=of J] {$\cdots$};
			\node (K2) [right=of K] {$\cdots$};
			\node (L2) [right=of L] {$\cdots$};
			\node (M2) [right=of M] {$\cdots$};
			\draw (A)--(B)--(C)--(D)--(H)--(H2);
			\draw (B)--(E)--(F)--(I)--(J)--(J2);
			\draw (F)--(G)--(K)--(L)--(M)--(M2);
			\draw (D)--(D2);
			\draw (I)--(I2);
			\draw (K)--(K2);
			\draw (L)--(L2);
		\end{tikzpicture}
		\caption{Initial structure of the directed rooted tree of congruence classes  $\bmod\, 2^{\sigma_n}$, up to $\mathbb{V}(4)$.}
		\label{fig:3}
	\end{figure}


	\section{Finite-Range Coverage and a Computable Bound}
	\label{sec:coverage}
	
	We examine the coverage of $\mathbb{Z}_{\ge 2}$ by the \emph{coefficient stopping time} classes generated up to order $N$. In addition to these classes, we include the trivial families: even integers satisfy $\sigma^\ast(x)=1$, and integers $x \equiv 1 \pmod 4$ satisfy $\sigma^\ast(x)=2$.
	
	\begin{defn}\label{def:coverage_set}
	    Fix $N\geq 1$. The \emph{coverage set of order $N$} is
	    \begin{align}\label{eq:coverage_set}
	        \mathcal{U}_N :=
	        & \{x\in\mathbb{Z}_{\geq 2}: x \equiv 0 \pmod 2\}\ \cup \nonumber
	        \\
	        & \{x\in\mathbb{Z}_{\geq 2}: x \equiv 1 \pmod 4\}\ \cup \nonumber
	        \\
	        & \bigcup_{n=1}^{N}\ \bigcup_{a\in\mathcal{C}_n}\ \{x\in\mathbb{Z}_{\geq 2}: x \equiv a \pmod {2^{\sigma_n}}\}.
	    \end{align}
	    The \emph{coverage bound} $K(N)$ is
	    \[
	        K(N):=\max\{K\geq 2 : \{2,3,\dots,K\}\subseteq \mathcal{U}_N\},
	    \]
	    where $K(N)=\infty$ if $\mathcal{U}_N=\mathbb{Z}_{\geq 2}$.
	\end{defn}
	
	\noindent Note that $\mathcal{U}_N$ refers exclusively to coefficient stopping times $\sigma^\ast(x)=\sigma_n$ with $n \le N$ (and trivial cases). This implies no statement about global dynamics.
	
	Lemma~\ref{lem:periodicity_UN} shows that membership in $\mathcal{U}_N$ is decidable via a finite check over one period $\bmod\, 2^{\sigma_N}$.
	
	\begin{lem}\label{lem:periodicity_UN}
	    Let $N\geq 1$ and $M:=2^{\sigma_N}$. The set $\mathcal{U}_N$ is periodic with period $M$, i.e.,
	    \[
	        x\in\mathcal{U}_N \iff x+M\in\mathcal{U}_N \qquad (x\in\mathbb{Z}_{\geq 2}).
	    \]
	\end{lem}
	
	\begin{proof}
		By Definition~\ref{def:coverage_set}, $\mathcal{U}_N$ is a union of residue classes with moduli $2$, $4$, and $2^{\sigma_n}$ for $1\le n\le N$. Since $\sigma_n\le \sigma_N$, each modulus divides $M=2^{\sigma_N}$. Therefore, if $x \equiv r \pmod m$ with $m\mid M$, then $x+M \equiv x \equiv r \pmod m$, and every class in \eqref{eq:coverage_set} is invariant under translation by $M$.
	\end{proof}
	
	\begin{cor}\label{cor:finite_check_UN}
	    Assume $\mathcal{U}_N\neq \mathbb{Z}_{\geq 2}$ and let $x_{\mathrm{miss}}:=\min\{x\geq 2: x\notin\mathcal{U}_N\}$. Then
	    \[
	        x_{\mathrm{miss}}\leq 2^{\sigma_N}+1 \qquad\text{and} \qquad K(N)=x_{\mathrm{miss}}-1.
	    \]
	    Testing membership in $\mathcal{U}_N$ reduces to the finite interval $[2, 2^{\sigma_N}+1]$.
	\end{cor}
	
	\begin{proof}
	    Set $M=2^{\sigma_N}$. By Lemma~\ref{lem:periodicity_UN}, membership in $\mathcal{U}_N$ depends only on the residue class $\bmod\, M$. If $x \notin \mathcal{U}_N$ exists, the set $\{x, x-M, \dots\}$ contains a representative in $\{2, \dots, M+1\}$ which is also missing. Thus $x_{\mathrm{miss}}\leq M+1$. The identity for $K(N)$ follows from Definition~\ref{def:coverage_set}.
	\end{proof}
	
	\begin{prop}[Exact coverage density]\label{prop:coverage_density}
	    Let $M_N=2^{\sigma_N}$. Then
	    \[
	        \#(\mathcal{U}_N\bmod M_N)
	        =3\cdot2^{\sigma_N-2}
	        +\sum_{n=1}^{N}|\mathcal{C}_n|2^{\sigma_N-\sigma_n}.
	    \]
	    Consequently, the density of $\mathcal{U}_N$ in one period is
	    \begin{equation}\label{eq:coverage_density}
	        \delta_N=\frac34+\sum_{n=1}^{N}\frac{|\mathcal{C}_n|}{2^{\sigma_n}}.
	    \end{equation}
	\end{prop}
	
	\begin{proof}
	    The even classes and the classes congruent to $1$ modulo $4$ are disjoint and occupy $M_N/2+M_N/4$ residues modulo $M_N$. By Theorem~\ref{thm_A100982_A076227}, the class families $\mathcal{C}_n$ describe distinct exact coefficient stopping times and are therefore disjoint. Each class modulo $2^{\sigma_n}$ lifts to $2^{\sigma_N-\sigma_n}$ classes modulo $M_N$. Summation gives the first formula, and division by $M_N$ gives \eqref{eq:coverage_density}.
	\end{proof}
	
	\begin{prop}[Meaning of the coverage bound]\label{prop:coverage_equivalence}
	    The sequence $K(N)$ is nondecreasing, and
	    \[
	        K(N)\longrightarrow\infty
	        \quad\Longleftrightarrow\quad
	        \sigma^\ast(x)<\infty\ \text{for every }x\geq2.
	    \]
	\end{prop}
	
	\begin{proof}
	    Since $\mathcal{U}_N\subseteq\mathcal{U}_{N+1}$, the sequence $K(N)$ is nondecreasing. If $K(N)\to\infty$, then every fixed $x\geq2$ belongs to some $\mathcal{U}_N$ and hence has finite coefficient stopping time.

	    Conversely, assume that every $x\geq2$ has finite coefficient stopping time. If $x$ is even or $x\equiv1\pmod4$, then $x$ belongs to the trivial families in $\mathcal{U}_1$. Otherwise $x>1$ is odd and $n(x)\geq1$. Lemma~\ref{lem:sigma_values} gives $\sigma^\ast(x)=\sigma_{n(x)}$, and Theorem~\ref{thm_A100982_A076227} yields
	    \[
	        x\bmod 2^{\sigma_{n(x)}}\in\mathcal{C}_{n(x)}.
	    \]
	    Thus every fixed $x$ belongs to $\mathcal{U}_N$ for all sufficiently large $N$. Applying this to the finite set $\{2,\ldots,K\}$ gives $K(N)\geq K$ eventually.
	\end{proof}
	
	\noindent Corollary~\ref{cor:finite_check_UN} allows computing $K(N)$ via finite verification. Appendix~\ref{sec:prog_table_5} implements this using the vectors $\mathbb{V}(n)$.
	
	Table~\ref{tab:coverage} shows that $x_{\mathrm{miss}}=27$ persists up to $N=18$. Direct iteration gives $\sigma^\ast(27)=59=\sigma_{36}$, so $27$ enters the coverage set only at order $36$. In particular, $27\notin \mathcal{U}_{18}$ despite a coverage density of $\approx 98.9\%$ in the period $M=2^{\sigma_{18}}$. Thus, high density does not imply coverage of an initial interval. Appendix~\ref{sec:prog_table_5} generates these values.
	
	\begin{table}[!ht]
	    \centering\small
	    \begin{tabular}{r r r r r r r r}
	        \hline
	        $N$ & $\sigma_N$ & $|\mathcal{C}_N|$ & $M$ & $\#\mathrm{cov}$ & $\#\mathrm{cov}/M$ & $x_{\mathrm{miss}}$ & $K(N)$\\
	        \hline
	        1 & 4 & 1 & 16 & 13 & 0.8125000000 & 7 & 6\\
	        2 & 5 & 2 & 32 & 28 & 0.8750000000 & 7 & 6\\
	        3 & 7 & 3 & 128 & 115 & 0.8984375000 & 27 & 26\\
	        4 & 8 & 7 & 256 & 237 & 0.9257812500 & 27 & 26\\
	        5 & 10 & 12 & 1024 & 960 & 0.9375000000 & 27 & 26\\
	        6 & 12 & 30 & 4096 & 3870 & 0.9448242188 & 27 & 26\\
	        7 & 13 & 85 & 8192 & 7825 & 0.9552001953 & 27 & 26\\
	        8 & 15 & 173 & 32768 & 31473 & 0.9604797363 & 27 & 26\\
	        9 & 16 & 476 & 65536 & 63422 & 0.9677429199 & 27 & 26\\
	        10 & 18 & 961 & 262144 & 254649 & 0.9714088440 & 27 & 26\\
	        11 & 20 & 2652 & 1048576 & 1021248 & 0.9739379883 & 27 & 26\\
	        12 & 21 & 8045 & 2097152 & 2050541 & 0.9777741432 & 27 & 26\\
	        13 & 23 & 17637 & 8388608 & 8219801 & 0.9798766375 & 27 & 26\\
	        14 & 24 & 51033 & 16777216 & 16490635 & 0.9829184413 & 27 & 26\\
	        15 & 26 & 108950 & 67108864 & 66071490 & 0.9845419228 & 27 & 26\\
	        16 & 27 & 312455 & 134217728 & 132455435 & 0.9868698940 & 27 & 26\\
	        17 & 29 & 663535 & 536870912 & 530485275 & 0.9881058242 & 27 & 26\\
	        18 & 31 & 1900470 & 2147483648 & 2123841570 & 0.9889907995 & 27 & 26\\
	    \end{tabular}
	    \caption{Finite-range coverage statistics within the period $M=2^{\sigma_N}$ for the set $\mathcal{U}_N$.}
	    \label{tab:coverage}
	\end{table}


	\section{Concluding Remarks}
	\label{sec:final_conclusion}
	
	The coefficient stopping-time classes admit two equivalent finite descriptions: by admissible positions of odd iterates and by residue classes modulo powers of $2$. The odd-position description gives a canonical parent map and hence a complete rooted-tree structure. It also yields the counting recursion without computational assumptions.
	
	The affine iterate formula provides explicit arithmetic relations between related classes. Within each fixed class, the coefficient stopping time and the classical stopping time are equal for all sufficiently large representatives, with the number of possible exceptions bounded explicitly in Corollary~\ref{cor:eventual_descent}. Any discrepancy between the two notions is therefore confined to a finite initial subset of each class.
	
	For every fixed order $N$, the union of the classes considered is periodic and has the exact density \eqref{eq:coverage_density}. The condition $K(N)\to\infty$ is equivalent to finite coefficient stopping time for every starting value. It would not, by itself, prove equality with the classical stopping time or exclude nontrivial cycles. These remain separate global questions.


	\appendix
	\section{Algorithms in PARI/GP}
	
	\textit{Note}: Program~\ref{sec:prog7} uses the global object \texttt{V} produced by Program~\ref{sec:prog3}. The remaining programs are standalone. Large values of \texttt{limit} may require increasing the PARI/GP stack size, for example with \texttt{default(parisize, "2G")}. The algorithms can also be run directly in a web browser using the GP/WASM interface\footnote{\url{https://pari.math.u-bordeaux.fr/gpwasm.html}}.
	
	\subsection{Algorithm for Proposition~\ref{prop:R_recurrence} and A100982}
	\label{sec:prog1}
	
	This program generates the data presented in Table~\ref{tab:1}, corresponding to sequence \OEIS{A100982}. For each $n\geq 2$, it outputs the entries of column $n$ along with their sum $|\mathcal{C}_n|$.
	
	\begin{lstlisting}
	{
		/* Define local variables */
		my(limit = 20);
		/* Calculate the maximum necessary length for the vectors based on log(3)/log(2). Adding a buffer (+5) to prevent out-of-bounds errors */
		my(max_rows = floor(limit * log(3)/log(2)) + 5);
		/* Index shift: v_prev[k+1] stores R(k,n-1). */
		my(v_prev = vector(max_rows));
		v_prev[2] = 1;                 /* R(1,1)=1 */
		for(n = 2, limit,
			my(Kappa = floor(n * log(3) / log(2)));
			my(Cn = 0);
			/* 'v_curr' represents the current column n. We re-initialize it in every iteration to ensure it starts with zeros. */
			my(v_curr = vector(max_rows));
			print();
			print1("For n=", n, " in column n: ");
			for(k = n, Kappa,
				/* Recurrence relation: R[k+1,n] = R[k,n] + R[k,n-1] */
				/* Implementation: New Value = (Value above in current col) + (Value in previous col) */
				v_curr[k+1] = v_curr[k] + v_prev[k];
				print1(v_curr[k+1], ", ");
				/* Update the sum using the increment operator */
				Cn += v_curr[k+1];
			);
			print(); 
			print(" and the sum is Cn=", Cn);
			/* IMPORTANT: Before moving to the next 'n', the current column becomes the "previous" one. */
			v_prev = v_curr;
		);
	}	\end{lstlisting}
	
	\subsection{Algorithm for Proposition~\ref{prop:R_recurrence} and A076227}
	\label{sec:prog2}
	
	This program generates the data presented in Table~\ref{tab:1}, corresponding to sequence \OEIS{A076227}. For each $k\geq 2$, it outputs the entries of row $k$ along with their sum $w(k)$.
	
	\begin{lstlisting}
	{
		/* Define local variables */
		my(limit = 20);
		/* Calculate the maximum necessary number of rows (buffer for Kappa) */
		/* Kappa grows faster than n, so we need more rows than columns. */
		my(max_rows = floor(limit * log(3)/log(2)) + 5);
		my(R = matrix(max_rows, limit));
		/* Index shift: R[k+1,n] stores R(k,n). */
		R[2, 1] = 1;                  /* R(1,1)=1 */
		for(n = 2, limit,
			/* Print header only for n > 2, matching the original logic */
			if(n > 2, 
				print(); 
				print1("For k=", n-1, " in row k: ");
			);
			my(Kappa = floor(n * log(3) / log(2)));
			/* Calculate values for column n */
			for(k = n, Kappa,
				R[k+1, n] = R[k, n] + R[k, n-1];
			);
			/* Calculate lower bound t (cf. OEIS A020915) */
			my(t = floor(1 + (n-1) * log(2) / log(3)));
			my(Wk = 0);
			/* Print row k=n-1 (accessing values from previous columns i). */
			if(n > 2,
				for(i = t, n-1,
					print1(R[n, i], ", ");
					Wk += R[n, i];
				);
				print(); 
				print(" and the sum is w(k)=", Wk);
			);
		);
	}	\end{lstlisting}
	
	\subsection{Algorithm for Theorem~\ref{thm_algo_parity_vectors}}
	\label{sec:prog3}

	This program generates the parity vectors of $\mathbb{V}(n)$ from the parent map in Theorem~\ref{thm_algo_parity_vectors}. The global object \texttt{V} is a vector of lists, with \texttt{V[n]} containing the vectors at level $n$ in generation order. Storage therefore grows with the actual number of vectors rather than with a fixed rectangular allocation.

	\begin{lstlisting}
	{
		/* V is global: V[n] is a List containing the vectors of V(n). */
		my(limit = 14);
		my(Log32 = log(3) / log(2));
		V = vector(limit);
		V[1] = List();
		listput(V[1], vector(2, i, 1));          /* root (1 1) */
		for(n = 2, limit,
			my(Kappa = floor(n * Log32));
			my(KappaPrev = floor((n-1) * Log32));
			V[n] = List();
			for(v = 1, #V[n-1],
				my(A = V[n-1][v]);
				my(B = vector(Kappa + 1));
				for(i = 1, KappaPrev + 1, B[i] = A[i]);
				/* Step 1: rightmost child, new 1 at position Kappa. */
				B[Kappa + 1] = 1;
				listput(V[n], B);
				/* Step 2: shift the new terminal 1 to the left. */
				for(j = 1, Kappa - n,
					if(B[Kappa + 1 - j] == 1, break);
					B[Kappa + 1 - j] = 1;
					B[Kappa + 2 - j] = 0;
					listput(V[n], B);
				);
			);
			print("n=", n, "   |C_n|=", #V[n]);
			my(p = 1);
			my(h0 = 2);
			for(i = 1, #V[n],
				my(B = V[n][i]);
				my(h = 0);
				while(h < #B && B[h+1] == 1, h++);
				if(h > h0, p = 1; h0 = h; print(); );
				print(B, "   ", h, "   ", p, "   ", i);
				p++;
			);
			print();
		);
	}	\end{lstlisting}

	\subsection{Algorithm for Table~\ref{tab:2} }
	\label{sec:prog4}
	
	This standalone program generates the counts $\mathcal{P}(h,n)$ in Table~\ref{tab:2}. It uses the same parent map and generation order as Program~\ref{sec:prog3}, but stores only one level at a time.
	
	\begin{lstlisting}
	{
	    my(Nmax = 11);
	    my(max_h = 12);
	    my(Log32 = log(3)/log(2));
	    my(Counts = matrix(max_h, Nmax));
	    Counts[2, 1] = 1;
	    my(Vprev = List());
	    listput(Vprev, [1, 1]);
	    for(n = 2, Nmax,
	        my(Kappa = floor(n*Log32));
	        my(Vcurr = List());
	        for(parent_index = 1, #Vprev,
	            my(parent = Vprev[parent_index]);
	            my(last_old = 0);
	            for(i = 1, #parent,
	                if(parent[i] == 1, last_old = i);
	            );
	            my(child = vector(Kappa + 1));
	            for(i = 1, #parent, child[i] = parent[i]);
	            child[Kappa + 1] = 1;
	            listput(Vcurr, child);
	            my(pos = Kappa + 1);
	            while(pos - 1 > last_old,
	                my(next_child = vector(Kappa + 1, i, child[i]));
	                next_child[pos - 1] = 1;
	                next_child[pos] = 0;
	                child = next_child;
	                pos--;
	                listput(Vcurr, child);
	            );
	        );
	        Vprev = Vcurr;
	        for(j = 1, #Vprev,
	            my(v = Vprev[j]);
	            my(h = 0);
	            while(h < #v && v[h + 1] == 1, h++);
	            if(h <= max_h, Counts[h, n]++);
	        );
	    );
	    for(h = 2, max_h,
	        print1("\\hline$\\TC{L}{", h, "}$");
	        for(n = 1, Nmax,
	            if(Counts[h, n] == 0,
	                print1("&");
	            ,
	                print1("&", Counts[h, n]);
	            );
	        );
	        print("&$\\cdots$\\\\");
	    );
	}	\end{lstlisting}
	
	\subsection{Algorithm for Lexicographically Ordered Permutations}
	\label{sec:prog5}
	
	This function generates the complete set of permutations of the binary word~\eqref{vector} in lexicographical order, corresponding to Table~\ref{tab:3}.
	
	\begin{lstlisting}
	NextPermutation(a) =
	{
	    my(i, k, t);
	    i = #a - 1;
	    while(i >= 1 && a[i] >= a[i+1], i--);
	    if(i < 1, return(0));
	    k = #a;
	    while(a[k] <= a[i], k--);
	    t = a[k]; a[k] = a[i]; a[i] = t;
	    /* reverse suffix */
	    for(k = i + 1, (#a + i)\2,
	        t = a[k];
	        a[k] = a[#a + 1 + i - k];
	        a[#a + 1 + i - k] = t;
	    );
	    return(a);
	};	\end{lstlisting}
	
	\subsection{Algorithm for Corollary~\ref{cor_xnp}}
	\label{sec:prog6}
	
	This program implements the rightmost-child transition in Corollary~\ref{cor_xnp} along the first branch of the generation order. It outputs the canonical representative $x$ together with the unique odd parameter $\lambda$.
	
	\begin{lstlisting} 
	{
		/* Define local variables */
		my(j = 3);
		my(xn = 3);
		my(limit = 20);
		my(Log32 = log(3)/log(2));
		my(A);
		/* Initial parity vector B */
		my(B = vector(j + 1));
		B[1] = 1; 
		B[2] = 1;
		for(n = 2, limit,
			my(Sigma = floor(1 + (n+1) * Log32));
			my(d = floor(n * Log32) - floor((n-1) * Log32));
			my(Kappa = floor(n * Log32));
			/* --- Step 1: Generate new parity vector B for n --- */
			if(n > 2,
				/* Create new vector B of size Kappa+1 */
				B = vector(Kappa + 1);
				/* Copy the prefix from the previous state A */
				for(i = 1, j - 1, B[i] = A[i]);
			);
			/* Apply logic based on d (modifying the 'frontier' at index j) */
			if(d == 2,
				B[j] = 0; 
				B[j+1] = 1;
				, /* else (d == 1) */
				B[j] = 1;
			);
			/* Update j for the next step */
			j += d;
			/* Save current state B to A for the next iteration (n+1) */
			A = B;
			/* --- Step 2: Determine indices Alpha where B[k] == 1 --- */
			my(Alpha = vector(n + 1));
			my(idx = 1);
			for(k = 1, Kappa + 1,
				if(B[k] == 1,
					Alpha[idx] = k - 1;
					idx++;
				);
			);
			/* --- Step 3: Solve Diophantine equation --- */
			my(r = Sigma - Kappa);          /* r in {2,3} */
			my(Lambda = 1, found = 0);
			/* Search for the unique valid odd Lambda */
			while(Lambda <= 2^r - 1,
			    my(x = xn + Lambda * 2^Kappa);
			    my(Sum_val = 0);
			    for(i = 1, n + 1, Sum_val += 3^(n+1-i) * 2^Alpha[i]);
			    my(y = (3^(n+1) * x + Sum_val) / 2^Sigma);
			    if(denominator(y) == 1,
			        x = lift(Mod(x, 2^Sigma));
			        print(n, " ", x, " ", Lambda);
			        xn = x;
			        found = 1;
			        break;
			    );
			    Lambda += 2;
			);
			if(!found, error("No valid lambda found at n = ", n));
		);
	}	\end{lstlisting}
	
	\subsection{Algorithm for Corollary~\ref{cor_delta_s}}
	\label{sec:prog7}

	This program validates Corollary~\ref{cor_delta_s} by detecting consecutive vectors in \texttt{V[n]} that differ by a single adjacent swap $01\leftrightarrow10$. It compares the predicted difference $\Delta S$ and the resulting relative congruence with direct computations. Transitions between different parent blocks are skipped. Program~\ref{sec:prog3} must be executed first with \texttt{limit} at least \texttt{Nmax}.

	\begin{lstlisting}
	/* S(v) = sum 3^(n+1-a) * 2^(alpha_a). */
	get_sv(v, n) = {
	    my(s = 0, a = 0);
	    for(i = 1, #v,
        	if(v[i] == 1, a++; s += 3^(n + 1 - a) * 2^(i - 1); );
    	);
    	return(s);
	};
	get_x_absolute(v, n, sigma) = {
	    my(m = 2^sigma);
	    my(inv3 = lift(Mod(3^(n+1), m)^-1));
	    return(lift(Mod(-inv3 * get_sv(v, n), m)));
	};
	detect_swap(v1, v2) = {
	    my(len = #v1);
	    if(#v2 != len, return(0));
	    my(diffs = List());
	    for(pos = 1, len, if(v1[pos] != v2[pos], listput(diffs, pos)); );
	    if(#diffs != 2, return(0));
	    my(i = diffs[1]);
	    if(diffs[2] != i + 1, return(0));
	    my(dir = 0);
	    if(v1[i] == 0 && v1[i+1] == 1 && v2[i] == 1 && v2[i+1] == 0, dir = -1);
	    if(v1[i] == 1 && v1[i+1] == 0 && v2[i] == 0 && v2[i+1] == 1, dir = +1);
	    if(dir == 0, return(0));
	    my(moved_pos = if(dir == -1, i+1, i));
	    my(k = 0);
	    for(pos = 1, moved_pos, if(v1[pos] == 1, k++));
	    return([i, dir, k]);
	};
	{
	    my(Nmax = 10);
	    if(type(V) != "t_VEC",
	        error("V not found. Run the parity-vector program first.");
	    );
	    my(Log32 = log(3)/log(2));
	    for(n = 2, Nmax,
	        my(sigma = floor(1 + (n + 1) * Log32));
	        my(M = 2^sigma);
	        my(inv3n1 = lift(Mod(3^(n+1), M)^-1));
	        my(tested = 0, ok = 0, fail = 0, skipped = 0);
	        for(i = 2, #V[n],
	            my(v_prev = V[n][i-1]);
	            my(v_curr = V[n][i]);
	            my(sw = detect_swap(v_prev, v_curr));
	            if(sw == 0, skipped++; next; );
	            tested++;
	            my(t = sw[1] - 1);
	            my(ds_expected = sw[2] * 3^(n + 1 - sw[3]) * 2^t);
	            my(ds_actual = get_sv(v_curr, n) - get_sv(v_prev, n));
	            my(x_prev = get_x_absolute(v_prev, n, sigma));
	            my(x_abs  = get_x_absolute(v_curr, n, sigma));
	            my(x_rel  = lift(Mod(x_prev - inv3n1 * ds_expected, M)));
	            my(cong_ok = 1);
	            if(t > 0,
	                if(lift(Mod(x_abs - x_prev, 2^t)) != 0, cong_ok = 0);
	            );
	            if(ds_expected != ds_actual || x_rel != x_abs || cong_ok == 0,
	                fail++;
	            ,
	                ok++;
	            );
	        );
	        print("n=", n, "   pairs=", #V[n]-1, "   tested=", tested,
	              "   OK=", ok, "   FAIL=", fail, "   SKIP=", skipped);
	    );
	}	\end{lstlisting}

	\subsection{Algorithm for Proposition~\ref{prop:terminal_recursion}}
	\label{sec:prog8}
	
	This program generates the canonical solutions $(x_n,y_n)$ on the terminal branch and the corresponding values $\beta_n$. In addition to the two recurrences, it verifies uniqueness in the canonical interval, the terminal Diophantine equation, and the direct residue formula modulo $2^{\sigma_n}$.
	
	\begin{lstlisting}
	{
	    my(limit = 20);
	    my(Log32 = log(3)/log(2));
	    my(xprev = 3, yprev = 2);
	    print(1, " ", xprev, " ", yprev, " --");
	    for(n = 2, limit,
	        my(Kappa = floor(n * Log32));
	        my(KappaNext = floor((n+1) * Log32));
	        my(d = KappaNext - Kappa);
	        my(sigma = KappaNext + 1);
	        my(M = 2^sigma);
	        my(candidates = 0, x, y, Beta);
	        for(b = 0, 5,
	            my(xb = (2*xprev - 1 + b*2^(Kappa+2))/3);
	            my(yb = (yprev + b*3^n)/2^(d-1));
	            if(denominator(xb) == 1 && denominator(yb) == 1
	                && xb >= 1 && xb < M,
	                candidates++;
	                x = xb;
	                y = yb;
	                Beta = b;
	            );
	        );
	        if(candidates != 1,
	            error("Canonical beta is not unique at n = ", n);
	        );
	        my(Sterminal = 3^(n+1) - 2^(n+1));
	        if(M*y != 3^(n+1)*x + Sterminal,
	            error("Terminal equation failed at n = ", n);
	        );
	        my(inv3 = lift(Mod(3^(n+1), M)^-1));
	        my(xabs = lift(Mod(-inv3*Sterminal, M)));
	        if(x != xabs,
	            error("Direct residue check failed at n = ", n);
	        );
	        print(n, " ", x, " ", y, " ", Beta);
	        xprev = x;
	        yprev = y;
	    );
	}	\end{lstlisting}
	
	\subsection{Algorithm for Table~\ref{tab:coverage}}
	\label{sec:prog_table_5}
	
	The coverage count is computed from Proposition~\ref{prop:coverage_density}; no vector of length $2^{\sigma_N}$ is allocated. The first missing value is found by testing the coefficient condition directly through time $\sigma_N$. Thus the program generates all rows of Table~\ref{tab:coverage}, including $N=18$, without constructing the class lists $\mathcal{C}_n$.
	
	\begin{lstlisting}
	/* Return the class counts |C_n| for 1 <= n <= Nmax. */
	class_counts(Nmax) = {
	    my(Log32 = log(3)/log(2));
	    my(max_rows = floor(Nmax*Log32) + 3);
	    my(R = matrix(max_rows, Nmax));
	    my(C = vector(Nmax));
	    /* R[k+1,n] stores R(k,n). */
	    R[2, 1] = 1;
	    C[1] = 1;
	    for(n = 2, Nmax,
	        my(Kappa = floor(n*Log32));
	        for(k = n, Kappa,
	            R[k+1, n] = R[k, n] + R[k, n-1];
	            C[n] += R[k+1, n];
	        );
	    );
	    return(C);
	};
	/* Test membership in U_N directly through time sigma_N. */
	covered_to_order(x, sigmaN) = {
	    my(z = x, m = 0);
	    for(k = 1, sigmaN,
	        if(z % 2,
	            m++;
	            z = (3*z + 1)/2;
	        ,
	            z = z/2;
	        );
	        if(2^k > 3^m, return(1));
	    );
	    return(0);
	};
	{
	    my(Nmin = 1);
	    my(Nmax = 18);
	    my(Log32 = log(3)/log(2));
	    my(C = class_counts(Nmax));
	    print("%% LaTeX rows for Table tab:coverage");
	    for(N = Nmin, Nmax,
	        my(sigmaN = floor(1 + (N+1)*Log32));
	        my(M = 2^sigmaN);
	        /* Exact count from the density formula. */
	        my(covCount = 3*M/4);
	        for(n = 1, N,
	            my(sigman = floor(1 + (n+1)*Log32));
	            covCount += C[n]*2^(sigmaN - sigman);
	        );
	        /* Search only until the first missing value is found. */
	        my(xmiss = 0);
	        for(x = 2, M + 1,
	            if(!covered_to_order(x, sigmaN),
	                xmiss = x;
	                break;
	            );
	        );
	        print1(N, " & ", sigmaN, " & ", C[N], " & ", M,
	               " & ", covCount, " & ");
	        printf("%.10f & ", covCount*1.0/M);
	        if(xmiss == 0,
	            print("-- & $\\infty$\\\\");
	        ,
	            print(xmiss, " & ", xmiss - 1, "\\\\");
	        );
	    );
	}	\end{lstlisting}


	\section*{Acknowledgments}
	The author is grateful for detailed comments on earlier versions that improved the clarity and correctness of the manuscript. AI-based tools were used for language editing, structural review, and computational and mathematical cross-checking. The author independently verified all statements and assumes full responsibility for the mathematical content.


\begin{thebibliography}{99}
		\bibitem{Everett77}
		Everett,~C.~J., \emph{Iteration of the number theoretic function $f(2n)=n$, $f(2n+1)=3n+2$}, Advances in Math. 25 (1977), 42--45.

		\bibitem{Garner81}
		Garner,~L.~E., \emph{On the Collatz 3n + 1 Algorithm}, Proc. Amer. Math. Soc., Vol. 82, No. 1 (May, 1981), pp. 19 -- 22, \url{https://www.jstor.org/stable/2044308}.
		
		\bibitem{Hikawa2026}
		Hikawa,~K., \emph{Finite-Dimensional Combinatorial and Arithmetic Structures of Parity Vectors for the Accelerated Collatz Map}, preprint, July 2026, \url{https://doi.org/10.13140/RG.2.2.29894.84804}.
		
		\bibitem{Lagarias85}
		Lagarias,~J.~C., \emph{The 3x + 1 problem and its generalizations}, The American Mathematical Monthly, Vol. 92, No. 1 (January, 1985), pp. 3--23, \url{https://www.cecm.sfu.ca/organics/papers/lagarias/paper/html/paper.html}.
		
		\bibitem{OEIS}
		Sloane,~N.~J.~A.~et al., \emph{The On-Line Encyclopedia of Integer Sequences},
		\url{https://oeis.org}.

		\bibitem{Terras76}
		Terras,~R., \emph{A stopping time problem on the positive integers}, Acta Arithmetica, Vol. 30 (1976), 241--252, \url{https://doi.org/10.4064/aa-30-3-241-252}.
		
		\bibitem{Terras78}
		Terras,~R., \emph{On the existence of a density} (Addendum), Acta Arithmetica, Vol. 35 (1978), 101--102.
		
		\bibitem{Thaler}
		Thaler,~M.~H., \emph{A $\delta$-recurrence for Winkler's Corollary 9}, 2025, \url{https://zenodo.org/records/17278225}.
		
		\bibitem{Winkler2018}
		Winkler,~M., \emph{The Recursive Stopping Time Structure of the 3x+1 Function}, March 2018 (uploaded in 2024), \url{https://arxiv.org/abs/1709.03385v4}.
		
		\bibitem{Wirsching98}
		Wirsching,~G.~J., \emph{The Dynamical System Generated by the $3n + 1$ Function}, Lecture Notes in Mathematics 1681, Springer-Verlag, Berlin Heidelberg, 1998, \url{https://link.springer.com/book/10.1007/BFb0095985}.
	\end{thebibliography}
\end{document}